\documentclass[11pt]{article}
\usepackage[colorlinks=true,
linkcolor=webgreen,
filecolor=webbrown,
citecolor=webgreen]{hyperref}

\usepackage{amssymb,psfig,epsfig,latexsym,graphicx,here}
\definecolor{webgreen}{rgb}{0,.5,0}
\definecolor{webbrown}{rgb}{.6,0,0}

\newtheorem{lemma}{Lemma}
\newtheorem{theorem}{Theorem}

\newtheorem{conj}{Conjecture}

\newcommand{\eqn}[1]{(\ref{#1})}
\newcommand{\bsq}{{\vrule height .9ex width .8ex depth -.1ex }}

\newcommand{\NN}{{\mathbb N}}

\newcommand{\sW}{{\mathcal W}}
\newcommand{\sS}{{\mathcal S}}
\newcommand{\sP}{{\mathcal P}}
\newcommand{\eeq}{\end{equation}}
\newcommand{\beql}[1]{\begin{equation}\label{#1}}

\makeatletter
\def\@sect#1#2#3#4#5#6[#7]#8{\ifnum #2>\c@secnumdepth
     \def\@svsec{}\else
     \refstepcounter{#1}\edef\@svsec{\csname the#1\endcsname.\hskip .75em }\fi
     \@tempskipa #5\relax
      \ifdim \@tempskipa>\z@
        \begingroup #6\relax
          \@hangfrom{\hskip #3\relax\@svsec}{\interlinepenalty \@M #8\par}%
        \endgroup
       \csname #1mark\endcsname{#7}\addcontentsline
         {toc}{#1}{\ifnum #2>\c@secnumdepth \else
                      \protect\numberline{\csname the#1\endcsname}\fi
                    #7}\else
        \def\@svsechd{#6\hskip #3\@svsec #8\csname #1mark\endcsname
                      {#7}\addcontentsline
                           {toc}{#1}{\ifnum #2>\c@secnumdepth \else
                             \protect\numberline{\csname the#1\endcsname}\fi
                       #7}}\fi
     \@xsect{#5}}
\def\@begintheorem#1#2{\it \trivlist \item[\hskip \labelsep{\bf #1\ #2.}]}

\def\section{\@startsection {section}{1}{\z@}{-3.5ex plus -1ex minus
 -.2ex}{2.3ex plus .2ex}{\normalsize\bf}}
\makeatother
\makeatletter
\def\subsection{\@startsection {subsection}{1}{\z@}{-3.5ex plus -1ex minus
 -.2ex}{2.3ex plus .2ex}{\normalsize\bf}}

\makeatother

\begin{document}
\begin{center}
{\large\bf The EKG Sequence} \\
\vspace*{+.5in}
{\em J. C. Lagarias, E. M. Rains} and {\em N. J. A. Sloane} \smallskip \\
Information Sciences Research Center \\
AT\&T Shannon Lab \\
Florham Park, NJ 07932--0971 \medskip \\
Email addresses: \href{mailto:jcl@research.att.com}{{\tt jcl@research.att.com}}, \href{mailto:rains@research.att.com}{{\tt rains@research.att.com}}, \href{mailto:njas@research.att.com}{{\tt njas@research.att.com}}
\bigskip \\
December 12, 2001; revised March 11, 2002 \bigskip

{\bf Abstract}
\end{center}

The EKG or electrocardiogram sequence is defined by $a(1) = 1$, $a(2) =2$ 
and, for $n \ge 3$, $a(n)$ is the smallest natural number 
not already in the sequence with the property that
\linebreak
${\rm gcd} \{a(n-1), a(n)\} > 1$.
In spite of its erratic local behavior, which
when plotted resembles an electrocardiogram,
its global behavior appears quite regular.  
We conjecture that almost all
$a(n)$ satisfy the asymptotic formula
$a(n) = n (1+ 1/(3 \log n)) + o(n/ \log n)$ as $n \to \infty$;
and that the exceptional values $a(n)=p$ and $a(n)= 3p$, for
$p$ a prime, produce the spikes in the EKG sequence.
We prove that $\{a(n): n \ge 1 \}$ 
is a permutation of the natural numbers and that
$c_1 n \le a (n) \le c_2 n$ for constants $c_1, c_2$. There remains
a large gap between what is conjectured and what is proved.

%
%

\section{Introduction}
\begin{figure}[htb]
\centerline{\includegraphics[angle=270, 
width=5in]{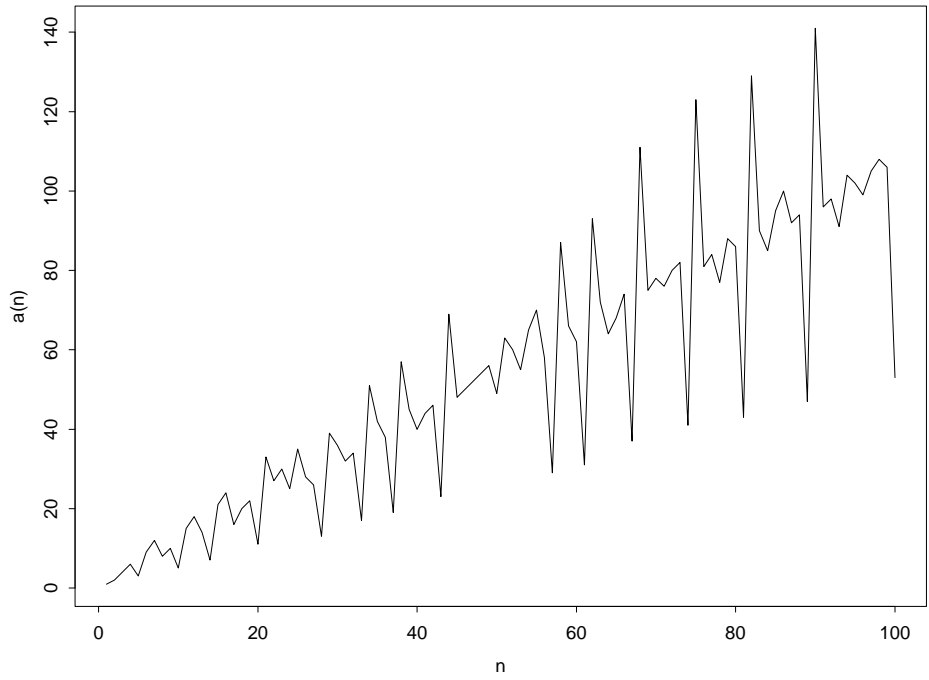}}

\caption{Plot of $a(1)$ to $a(100)$, with successive points joined by lines.}
\label{F1}
\end{figure}

\begin{figure}[H]
\centerline{\includegraphics[angle=270,
width=5in]{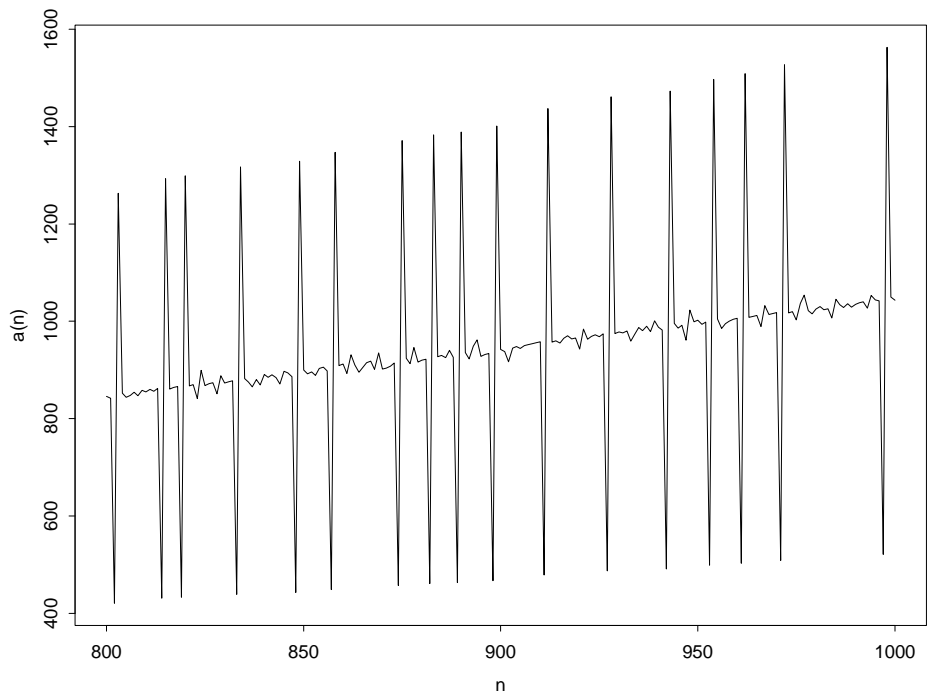}}

\caption{Terms 800 to 1000, with successive points joined by lines.}
\label{F2}
\end{figure}

Consider the sequence defined by
$a(1) =1$, $a(2) =2$ and, for $n \ge 3$, $a(n)$ is the smallest 
natural number not in $\{a(k) : 1 \le k \le n-1\}$ with the 
property that ${\rm gcd} \{a(n-1), a(n) \} \ge 2$.
This sequence might be called a greedy gcd sequence, but
because of its striking appearance when plotted we will name it the
{\em EKG} (or {\em electrocardiogram}) {\em sequence}---see
Figures \ref{F1}, \ref{F2}.
It was apparently first discovered by Jonathan Ayres
[Ayres 2001] and appears as sequence \htmladdnormallink{A064413}{http://www.research.att.com/cgi-bin/access.cgi/as/njas/sequences/eisA.cgi?Anum=A064413} in [Sloane 2001].
The first 30 terms are
$$
\begin{array}{rrrrrrrrrrr}
1 & 2 & 4 & 6 & 3 & 9 & 12 & 8 & 10 & 5 & ~ \\
15 & 18 & 14 & 7 & 21 & 24 & 16 & 20 & 22 & 11 & ~ \\
33 & 27 & 30 & 25 & 35 & 28 & 26 & 13 & 39 & 36 & \cdots
\end{array}
$$
Although the local behavior is erratic, 
plots of the first 1000 or 10000 terms show 
considerable regularity (see Figures \ref{F3}, \ref{F4} in Section 4).

The EKG sequence has a simple recursive definition, yet seems
surprisingly difficult to analyze.
Its definition
combines both additive and multiplicative aspects of the integers,
and the greedy property of its definition produces a complicated
dependence on the earlier terms of the sequence.
Indeed, it is not immediately obvious whether it contains all 
positive integers, but we show this is the case---the
EKG sequence  is a permutation of the positive integers.

In comparing Figures~\ref{F1} and \ref{F3}, one
is reminded of the contrast between the irregular plot
of $\pi (x)$ (the number of primes $\le x$) for $x \le 100$ 
and the very smooth plot for $x \le 50000$ as shown in Don Zagier's 
lecture on ``The first 50 million prime numbers'' [Zagier 1977].
This is not a  coincidence, as we will see, because (experimentally)
the spikes in the EKG sequence are associated with the primes arranged
in increasing order. However, the spacings between spikes are not the
same as the spacings between consecutive primes.

Although the EKG sequence itself seems only to have been 
proposed recently, in the early 1980's 
Erd\H{o}s, Freud and Hegyvari [Erd\H{o}s et~al. 1983]
studied properties of integer permutations with 
restrictions placed on allowed values of greatest common divisors of 
consecutive terms.

In Section 2 we derive a number of basic properties of the EKG sequence, 
and prove that it is a permutation of the natural numbers $\NN$.
An efficient algorithm for computing the sequence is given in Section 3.
Using this algorithm we computed $10^7$ terms;
this led us to conjectured asymptotic formulae given in Section 4.
We give a heuristic argument why these
formulae may be true, but it seems likely they will
be very hard to prove.  We  are able to rigorously  extablish
linear upper and lower bounds on the EKG sequence,
namely $\frac{1}{260} n \le a(n) \le 14n$
(Sections 5 and 6); the proofs use sieving ideas. 
Section 7 discusses experimental results
concerning the cycle structure of the associated permutation 
of $\NN$. The final section discusses generalizations to
other sorts of integer permutations resulting from
greedy constructions with restrictions on gcd's of consecutive
terms.

%
%

\section{The sequence is a permutation}

We begin with some general remarks about the sequence.

For $n \ge 2$, let $g = {\rm gcd} \{a(n-1), a(n) \}$.
For some prime $p$ dividing $a(n-1)$, $a(n)$ is the 
smallest multiple of $p$ not yet seen (otherwise the smaller
multiple of $p$ would be a better candidate for $a(n)$).
We call such primes $p$ the {\em controlling primes} for $a(n)$.
There may be more than one, and their product divides $g$.

For any prime $p$ and $n \ge 2$, let $B_p (n)$ be the smallest multiple of 
$p$ that is not in
\linebreak
$\{a(1), \ldots, a(n-1) \}$.
For example, the sequence $\{B_2 (n): n \ge 2 \}$ 
begins 2, 4, 6, 8, 8, 8, 8, 10, 14, $\ldots$.
Clearly
\beql{EqB1}
B_p (n) \le B_p (n+1) \le pn
\eeq
for all $p,n \ge 2$.
Then we have $a(1) =1$,
\beql{Eq1}
a(n) = \min \{B_p (n): ~\mbox{$p$ divides $a(n-1)$} \} ~,
\eeq
for $n \ge 2$, which provides an alternative definition of the sequence.
\begin{lemma}\label{L1}
Let $p$ be a prime $>2$ that divides some term of the sequence.
If $p$ first divides $a(n)$ then $a(n) = qp$ where $q$ is the 
smallest prime dividing $a(n-1)$, $q$ is less than $p$,
$a(n+1) =p$, and either $a(n)$ or $a(n+2)$ is equal to $2p$.
The new primes that divide the terms of the sequence appear 
in increasing order.
\end{lemma}

\paragraph{Proof.}
Let $a(n)$ be the first term divisible by $p$.
The numbers $pq$ where $q$ is a prime dividing $a(n-1)$ are 
all candidates for $a(n)$, and so $a(n) = pq$ where $q$ 
is the smallest such prime.
Also $p$ must be the smallest prime that has not appeared as 
a divisor of $\{a(1), \ldots, a(n-1) \}$ (for if $p'$ were a 
smaller such prime then $p'q$ would be a better candidate for $a(n)$).
In particular the primes
that divide the terms of the sequence must appear in increasing order, 
and $q < p$.
Then $p$ is a candidate for $a(n+1)$, and is less than 
$B_q (n+1) \ge B_q (n) = pq$, so $a(n+1) =p$.
Finally, either $a(n) =2p$ or else $2p$ is the winning candidate 
for $a(n+2)$.~~~$\bsq$

\begin{lemma}\label{L1a}
The primes that appear in the sequence occur in increasing order.
\end{lemma}

\paragraph{Proof.}
This follows from Lemma \ref{L1}, since the first time $p$ divides 
a term of the sequence the next term is $p$ itself.~~~$\bsq$

\begin{lemma}\label{L2}
If infinitely many multiples of a prime $p$ appear in the sequence 
then all multiples of $p$ appear.
\end{lemma}

\paragraph{Proof.}
We argue by contradiction, and let $kp$ be the first multiple of $p$ 
that is missed.
Choose $n_0$ so that $a(n) > kp$ for all $n \ge n_0$.
Since infinitely many multiples of $p$ occur, there exists
$n > n_0$ with $a(n) = lp$ for some $l$.
But now we must have $a(n+1) = kp$, because 
${\rm gcd} \{a(n), kp\} \ge p$ is allowed, and all 
smaller possible values which are ever going to appear in the 
sequence have already appeared.
This is a contradiction.~~$\bsq$

\begin{lemma}\label{L3}
If all multiples of a prime $p$ appear in the sequence then all
 positive integers appear.
\end{lemma}

\paragraph{Proof.}
Again we argue by contradiction and let $k \ge 2$ be the first integer 
that is missed.
Since infinitely many multiples of $k$ occur among all the 
multiples of $p$, we get a contradiction just as in Lemma \ref{L2}.
Namely, there exists for $n > n_0$ a value $a(n) = klp$ for some $l$, 
and ${\rm gcd} \{a(n),k \} \ge k$ is allowed, and all 
smaller possible values have already been used.
Thus $a(n+1) =k$, a contradiction.~~~$\bsq$

\begin{theorem}\label{T1}
$\{a(n) : n \ge 0\}$ is a permutation of the natural numbers.
\end{theorem}

\paragraph{Proof.}
No number can appear twice, by construction, so it suffices to show 
that every number appears.
Suppose only finitely many different primes divide 
the terms of the sequence. Then one of them would appear infinitely many times,
and Lemmas \ref{L2} and \ref{L3} would imply that all integers occur, which
is a contradiction.

Therefore infinitely many different primes $p$ divide the terms of the sequence.
Then by Lemma \ref{L1} infinitely many even numbers $2p$ occur,
by Lemma \ref{L2} all even numbers occur and by Lemma \ref{L3} 
all positive integers occur.~~~$\bsq$

\paragraph{Remark.} As will be discussed in Section 8, the
principle of this
 proof generalizes to a wide variety of other integer sequences
 defined by restrictions on the gcd's of consecutive terms.

%
%

\section{Numerical investigations}
To compute the EKG sequence it is better not to use the original definition
 but to use \eqn{Eq1} and to store the current values of $B_p (n)$ 
for primes $p$.
An efficient way to arrange the computation is to maintain four tables:
\begin{center}
\begin{tabular}{lll}
hit$(m)$ & = & 0 if $m$ has not yet appeared, otherwise 1; \\ [+.1in]
gap$(m)$ & = & current value of $B_m (n)$ if $m$ is a prime, 
otherwise $m$; \\ [+.1in]
small$(m)$ & = & smallest prime factor of $m$; \\ [+.1in]
quot$(m)$ & = & largest factor of $m$ not divisible by small$(m)$.
\end{tabular}
\end{center}
Combining these tables in a $C$ ``struct'' minimizes memory access.

Suppose we wish to compute the sequence until $a(n)$ reaches or exceeds $N$.
The first step is to precompute small$(m)$ and quot$(m)$ for $m \le N$.
Since it is only necessary to consider primes $\le \sqrt{N}$,
this takes about $\Sigma_{p \le \sqrt{N}} N/p = O(N \log \log N)$ steps.

In the main loop, let $a(n)$ be the current value.
Set $k = a(n)$, $B=N$ and repeat until $k$ reaches 1:
$$
\begin{array}{lll}
p & = & \mbox{small} (k) \,, \\
B & = & \min \{B, \mbox{gap} (p) \} \,, \\
k & = & \mbox{quot} (k) \,.
\end{array}
$$
Then we set $a(n+1) =B$, $\mbox{hit} (B) =1$, and update 
$\mbox{gap} (q)$ for primes $q$ dividing $B$.

We have not analyzed the complexity of the main loop in detail, 
but it also appears to take roughly $O(N \log \log N )$ steps, 
comparable to and not much greater than the number of steps 
needed for the precomputation part of the calculation.
The program computed $10^7$ terms of the sequence in less than a minute.
For example $a(10954982) = 11184814$.

%
%

\section{A conjectured asymptotic formula}

The results from the experimental data suggest that whenever a prime $p$ 
occurs in the sequence it is preceded by $2p$
and (consequently) followed by $3p$.
Although it is theoretically possible that some other 
multiple of $p$ occurs before $p$, for example, 
we might have seen $\ldots$, $3p$, $p$, $2p, \ldots$, 
this does not happen in the first $10^7$ terms.
\begin{conj}
\label{Con0}
Whenever a prime $p$ occurs in the sequence it is immediately 
preceded by $2p$ ${\rm (}$and hence followed by $3p{\rm )}$.
\end{conj}

\begin{figure}[htb]
\centerline{\includegraphics[angle=270, 
width=5in]{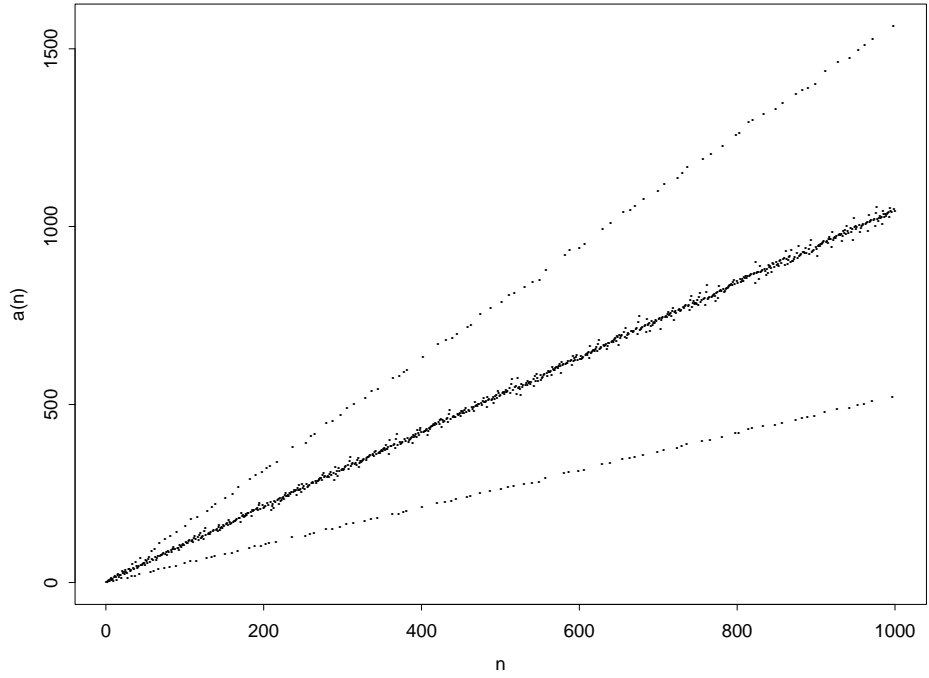}}

\caption{The first 1000 terms (represented by dots), successive points not joined.}
\label{F3}
\end{figure}

\begin{figure}[H]
\centerline{\includegraphics[angle=270, width=5in]
{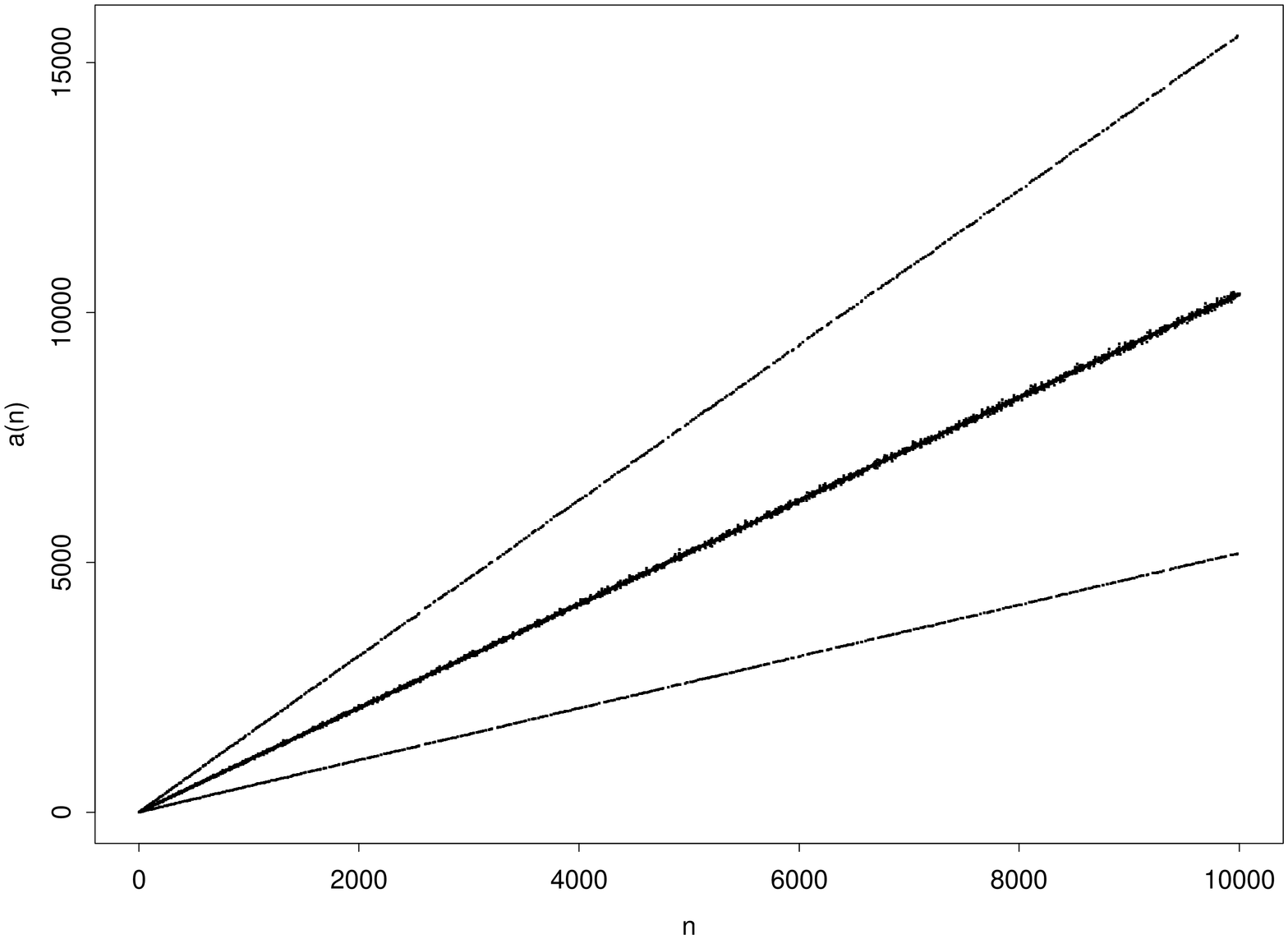}}

\caption{The first 10000 terms (represented by dots), successive points not joined.}
\label{F4}
\end{figure}

The numerical results also strongly suggest that the terms 
of the sequence fall close to three lines (see Figs. \ref{F3}, \ref{F4}).

\begin{figure}[htb]
\centerline{\includegraphics[angle=270, 
width=5in]{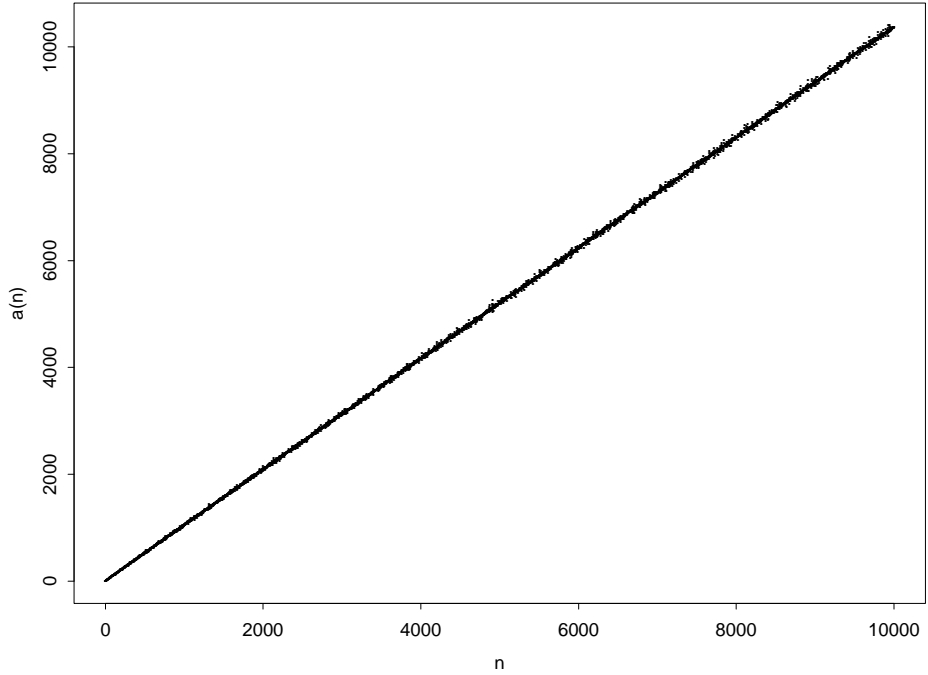}}

\caption{The sequence smoothed by replacing $a(n)=p$ or $3p$, 
$p$ prime $> 2$, by $a(n) =2p$.}
\label{F5}
\end{figure}

\begin{itemize}
\item
if $a(n) =m$ and $m$ is neither a prime nor three times a prime, 
then $a(n) \approx n$;
\item
if $a(n) =p$, $p$ prime, then $a(n) \approx n/2$;
\item
if $a(n) =3p$, $p$ prime, then $a(n) \approx 3n/2$.
\end{itemize}
This was also observed by Ayres [Ayres 2001].
In fact, if we smooth the sequence by replacing every term 
$a(n) = p$ or $3p$, $p$ prime $> 2$, by $a(n) =2p$,
the terms of the sequence lie close to a single line (see Fig. \ref{F5}).

\begin{figure}[htb]
\centerline{\includegraphics[angle=270, 
width=5in]{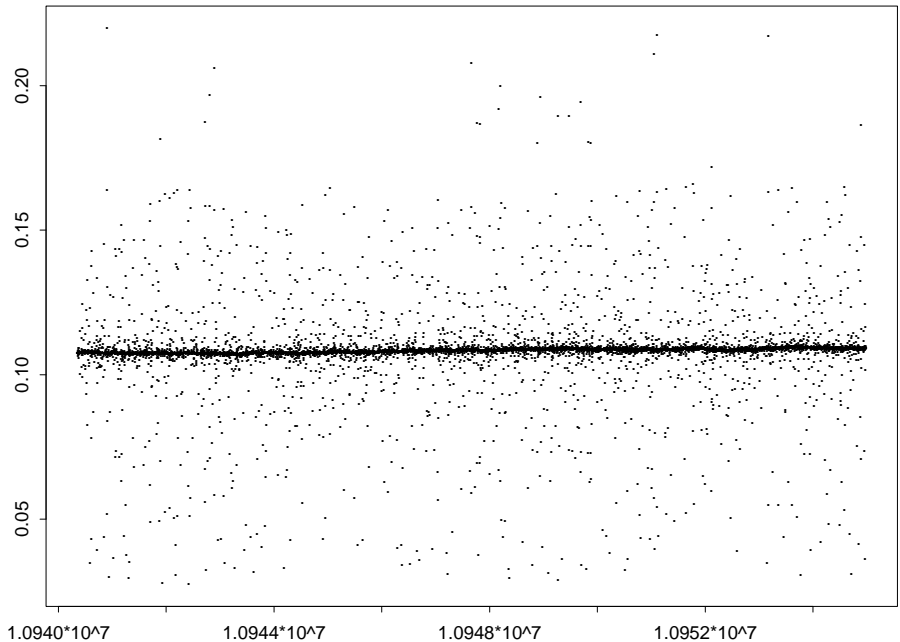}}

\caption{Values of $c$ in \eqn{Eq6} for $n$ near $10^7$.}
\label{F6}
\end{figure}

A plausible but nonrigorous argument suggests a more precise conjecture.

\begin{conj}\label{Con1}
Let $f(n) \sim g(n) $ mean that the ratio of the two sides 
approaches ${\rm 1}$ as $n \to \infty$.

${\rm (1)}$ If $a(n) = m$, $m \neq p$ or $3p$ for $p$ prime, then
\beql{Eq2}
a(n) \sim n \left( 1 + \frac{1}{3 \log n} \right) \,;
\eeq

${\rm (2)}$ If $a(n) =p$, $p$ prime, then
\beql{Eq3}
a(n) \sim \frac{1}{2} n \left( 1 + \frac{1}{3 \log n} \right) \,;
\eeq

${\rm (3)}$ If $a(n) = 3p$, $p$ prime, then
\beql{Eq4}
a(n) \sim \frac{3}{2} n \left( 1+ \frac{1}{3 \log n} \right) \,.
\eeq
\end{conj}

To see why this conjecture might be true,
consider a term $a(n) =m$ of the smoothed sequence.
Examination of Fig. \ref{F5} suggests that the smoothed sequence 
has hit all the numbers from 1 to $m$ at least once 
(the numbers occur a little out of order, but never by much).
However, we have smoothed away the numbers $p$ and $3p$ that are 
$\le m$, while picking up the primes $p$ and $3p$ for $2p \le m$.
Therefore
\beql{Eq5}
n \sim m - \pi (m) - \pi \left( \frac{m}{3} \right) + 2 \pi \left( \frac{m}{2} \right) \,,
\eeq
where $\pi (x) =$ number of primes $\le x$.
Then \eqn{Eq1} follows at once from the asymptotic formula 
$\pi (x) \sim x/ \log x$.
Equations \eqn{Eq3} and \eqn{Eq4} are based on \eqn{Eq1} and the 
observations made at the beginning of this section.

Although we are unable to prove this conjecture, 
it is an excellent fit to the data.

If we try to write
\beql{Eq6}
a(n) \approx n \left( 1 + \frac{1}{3 \log n} +
\frac{c}{(\log n)^2} \right) \, ~~~\mbox{(?)}
\eeq
then the values of $c$ do not appear to converge to a single value (see
Fig. \ref{F6}), although $c$ is very often close to $0.11$.
It seems conceivable that  $c$ might converge in distribution to a limiting
distribution.
Using two terms of the asymptotic expansion of $\pi (x)$ would give
\beql{Eq3b}
a(n) \sim n \left( 1 + \frac{1}{3 \log n} + \frac{c'}{(\log n)^2} \right) \, ~~~\mbox{(?)}
\eeq
where $c' = 4/9 + (\log 3)/3 - \log 2 = 0.1175\ldots$ .

Conjectures \ref{Con0} and \ref{Con1} predict that the 
$k$-th prime $p_k$ will occur in the pattern $a(n) = 2p_k$,
$a(n+1) = p_k$, $a(n+2) = 3p_k$, where
$$n \sim \frac{2p_k}{1+ \frac{1}{3 \log (2p_k )}}~.$$

These conjectures may be hard to settle, because
the permutation $a(n)$ encodes an intricate interaction between additive  
and multiplicative properties of integers, which by the
``greedy'' property of the definition depends on
all the earlier terms of the sequence. 

In the next two sections
we establish linear upper and lower bounds for the sequence,
namely 
$$\frac{1}{260}n \le a(n) \le 14n.$$
The numerical evidence supports the
following conjectural bounds.

\begin{conj}\label{Con3}
The sequence $a(n)$ satisfies
$a(n) \ge \frac{13}{28} n$, with equality if and only if $n=28$,
and $a(n) \le \frac{12}{7}n$, with equality if and only if $n=7$.
\end{conj}

For large $n$ the asymptotic lower and upper bounds
would be $n/2$ and $3n/2$, by  Conjecture \ref{Con1}.

%
%

\section{A linear upper bound}

In this section we show:

\begin{theorem}\label{T2}
$a(n) \le 14n$, for $n \ge 1$.
\end{theorem}

\paragraph{Remark.} The proof is by contradiction,
and the basic idea of the proof uses the fact that
the function values $\{a(k): 1 \le k \le n\}$ 
have a sandpile-like structure, in 
which each element can be built only at the top of  a ladder of all
smaller multiples of a controlling prime $p$ dividing it. 
To reach a number at height exceeding $14n$ via a ladder of
multiples of a controlling  prime $p$, 
one repeatedly falls off this ladder
while building it, at various smaller multiples $kp$ of $p$ of
size between $2n$ and $14n$ (see property (P4) below).
To get back on this ladder one must use ladders of other
controlling primes $q$ which reach to some currently omitted 
multiple of $p$ at the
top of their ladder. The total number of elements in such ladders is
shown to be large using a combinatorial sieve 
argument (compare [Hooley 1976], pp. 4-5);
a contradiction results by showing that the sandpile contains more than
$n$ elements. \\

We begin with a preliminary lemma.

\begin{lemma}\label{L5A}
If $a(n)$ is divisible by a prime $p$, then $p \le n$.
If $p \neq 2$, $p < n$.
\end{lemma}

\paragraph{Proof.}
The result is true if $p=2$ or $n \le 3$,
so we may assume $p \ge 3$ and $n \ge 4$.
We may also assume $a(n)$ is the first term divisible
by $p$.
By Lemma \ref{L1}, $a(n-1) = pq$ where $q$ is the 
controlling prime for $a(n)$.
Then $pq = a(n) = B_q (n) \le q(n-1)$ by \eqn{EqB1}.~~~$\bsq$

\paragraph{Proof of Theorem ~\ref{T2}.}
We argue by contradiction, and let $n$ be the smallest number 
such that $a(n) > 14n$.
By direct verification we know $n$ is large (in fact $n > 10^7$, 
although we will not use that in the proof).
Let $p$ be the smallest controlling prime for $a(n)$, say $a(n) = lp > 14n$,
so $l \ge 15$.
Also $a(n) = B_p (n)$, so from
\eqn{EqB1} and Lemma \ref{L5A}, $17 \le p < n$.

Let $l' = \lceil l/2 \rceil \ge 8$, and consider the ``window''
$$\sW = \{ l' p , (l'+1 ) p, \ldots, (l-1) p \} \,.
$$
Since $p$ is a controlling prime for $a(n)$, 
every element of $\sW$ has already appeared in the sequence:
$a(t) \in \sW$ implies $t \le n-1$.
Call $a(t) \in \sW$ an {\em entry point} if $a(t-1) \not\in \sW$, 
an {\em exit point} if $a(t+1) \not\in \sW$.
The following properties hold.

\noindent{(P1)}
The number of entry points is equal to the number of exit points.

\noindent{(P2)}
At most one entry point $a(t)$ has $p$ as a controlling prime.
This can only happen if
$a(t-1) = (l'-1) p$, $a(t) =l'p$, for some $t \le n-1$.

\noindent{(P3)}
At most one exit point $a(t)$ has $p$ as a controlling prime for $a(t+1)$.
This happens just if $a(n-1) \in \sW$.
Furthermore, if this is not the case then $a(n-1) = ip$ with
$i < l'$, and hence there is no entry point with $p$ as a controlling prime.

\noindent
{(P4)}
Let $a(t) = \alpha p \in \sW$,
$a \in [l', l-1 ]$.
If $\alpha$ is a multiple of a prime $q \le 7$, then $a(t)$ is an exit point.
(For $a(t+1) \le B_q(t+1) \le 7t \le 7(n-1) < l'p$.)

For a set of primes $\sP$, let $D_{\sP} (a,b)$ denote the number of 
integers $a \le i \le b$ such that $i$ is a multiple of some element of $\sP$.
Let
$$\theta = D_{\{ 2,3,5,7\}} (l',l-1)\,.$$
By (P4), the number of exit points is at least $\theta$, while
if there exists an entry point controlled by $p$ then there are in fact
at least $\theta +1$ exit points.  By (P1) we conclude
that the number of entry points not 
controlled by $p$ is at least $\theta$.

Let $\sS = \{q_1 , q_2, \ldots \}$ (with $q_1 < q_2 < \cdots$) denote the set
of controlling primes for entry points not controlled by $p$.
Note that $2,3,5,7,p$ are not in $\sS$ (by the same argument as in (P4)),
and $|\sS | \le \pi (l-1) -4$.
Then
\begin{eqnarray}\label{Eq8a}
\theta & = & D_{\{2,3,5,7\}} (l', l-1) \nonumber \\
& \le & \mbox{number of entry points not controlled by $p$} \nonumber \\
& \le & D_{\sS} (l', l-1 ) \nonumber \\
& \le & \sum_{q \in \sS} \left\lceil \frac{l-l'}{q} \right\rceil \nonumber \\
& \le & (l-l'-1) \sum_{q \in \sS} \frac{1}{q} + | \sS | \,.
\end{eqnarray}
Setting
$\phi = \sum_{q \in \sS} \frac{1}{q}$, we have
\beql{Eq8}
\phi \ge \frac{\theta - 
\pi (l-1) +4}{\left\lfloor \frac{l}{2} \right\rfloor -1} \,.
\eeq
The right-hand side of \eqn{Eq8} is a function of the single variable $l$,
and is $\ge 2/9$ for all $l \ge 15$, with equality 
if and only if $l$ is 20 or 21.
(This is easily verified by computer for small $l$, 
say $l \le 1000$, and analytically for larger $l$.)
In other words,
$$\sum_{i \ge 1} \frac{1}{q_i} \ge \frac{2}{9} \,.$$
Define $k$ by
\beql{Eq8b}
\sum_{i=1}^{k-1} \frac{1}{q_i} < \frac{2}{9} \le \sum_{i=1}^k \frac{1}{q_i} \,,
\eeq
and let $\sS' = \{ q_1, \ldots, q_k \} \subseteq \sS$.
Note that $q_k \ge 17$, since $1/11 + 1/13 + 1/17 \ge 2/9$, 
but no proper subset of $\{ 1/11$, $1/13$, $1/17 \}$ has this property.

Every multiple of any element $q \in \sS$ that is $\le l' p$ 
must have already occurred in the sequence, by definition.
We obtain a contradiction by showing that there are more than $n$ different
multiples of elements of $\sS$ that are $\le l' p$.
By inclusion-exclusion we have
\beql{Eq9}
D_\sS (1, l' p ) \ge D_{\sS'} (1, l'p ) \ge
\sum_{i=1}^k \left\lfloor \frac{l' p}{q_i} \right\rfloor -
\sum_{1 \le i < j \le k} \left\lfloor \frac{l' p}{q_i q_j} \right\rfloor \,.
\eeq
To bound the first term in \eqn{Eq9}, observe that for 
$q_i \in \sS$ there is a multiple of $q_i p$ in $\sW$, so
$$q_i < \frac{14n}{17} < \frac{lp}{17} \,,$$
\begin{eqnarray*}
\left\lfloor \frac{l' p}{q_i} \right\rfloor & \ge & 8 \,, \\
\left\lfloor \frac{l'p}{q_i} \right\rfloor
& \ge & \frac{8}{9} \frac{l'p}{q_i} \,,
\\
\sum_{i=1}^k \left\lfloor \frac{l'p}{q_i} \right\rfloor & \ge &
\frac{8}{9} \sum_{i=1}^k \frac{l'p}{q_i} \,.
\end{eqnarray*}
To bound the second term in \eqn{Eq9} we use
\begin{eqnarray*}
\lefteqn{\sum_{i\le i < j \le k} \frac{1}{q_i q_j}} \\
&& = \frac{1}{2}\left( \sum_{i=1}^{k-1} \frac{1}{q_i} \right)^2 - \frac{1}{2}
\sum_{i=1}^{k-1} \frac{1}{q_i^2} + \frac{1}{q_k}
\sum_{i=1}^{k-1} \frac{1}{q_i} \\
&& < \frac{1}{2} \left( \sum_{i=1}^{k-1} \frac{1}{q_i} \right)^2 +
\frac{1}{2q_k} \sum_{i=1}^{k-1} \frac{1}{q_i} \,.
\end{eqnarray*}
Then, since $q_k \ge 17$, we have
\begin{eqnarray*}
\lefteqn{D_\sS (1, l'p ) > l'p \left\{
\frac{8}{9} \left( \sum_{i=1}^k \frac{1}{q_i} \right) \right.} \\
&& \left. - \frac{1}{2}\left( \sum_{i=1}^{k-1} \frac{1}{q_i} \right)^2 -
\frac{1}{2q_k} \sum_{i=1}^{k-1} \frac{1}{q_i} \right\} \\
&& \ge l' p \left\{
\frac{8}{9} \cdot \frac{2}{9} - \frac{1}{2} \left( \frac{2}{9} \right)^2 -
\frac{1}{34} \cdot \frac{2}{9} \right\} \\
&& > 7n \frac{229}{1377} = \frac{1603}{1377} n > n \,,
\end{eqnarray*}
which is the desired contradiction.~~~$\bsq$

%
%

\section{A linear lower bound}

In this section we show:

\begin{theorem}\label{T3}
$a(n) \ge \lceil \frac{1}{260} n \rceil$, for $n \ge 1$.
\end{theorem}

\paragraph{Remark.}
The proof is a modification of that of Theorem \ref{T1}.
It aims to show that 
if some number less than $n/260$ is missed in
$\{ a(k): 1 \le k \le n\}$ then
there are at least $n/65$ numbers in this set that
are even numbers, and Lemma~\ref{L5} below provides the mechanism
to get a contradiction.
The method of Theorem \ref{T1} seems inherently weaker when used for a 
lower bound, so we have not attempted to streamline the proof.
It would certainly be possible to reduce the constant 260, 
but not to anything close to 14. \\

We begin with three lemmas.

\begin{lemma}\label{L4}
For a prime $p$, if $a(n) = kp$ for some $k$ then $a(j) = k$ 
for some $j \le n+1$.
\end{lemma}

\paragraph{Proof.}
We argue by induction on $k$.
The result is true for $k=1$ since $a(1) =1$.

Let $q$ be a controlling prime for $a(n)$, so $q | a(n-1)$ and $q| a(n) = kp$.
Case (i): $q \neq p$.
Then $q|k$, say $k= mq$, $a(n) = mqp$.
Hence all multiples $iq$ with $i < mp$ have already appeared, 
and in particular
$a(j) = mq =k$ for some $j < n$.

Case (ii):
$q=p$.
All multiples $ip$ with $i < k$ have already occurred.
By the induction hypothesis, $i$ has occurred for all $i < k$.
If $k$ has occurred then $a(j) =k$ with $j < n$, and otherwise
$a(n+1) =k$.~~~$\bsq$

\begin{lemma}\label{L5}
If at least $4k$ even numbers occur in $\{a(1),\ldots, a(n) \}$ 
then all numbers $\{1, \ldots, k \}$ occur in $\{a(1),\ldots, a(n+1) \}$.
\end{lemma}

\paragraph{Proof.}
In view of Lemma \ref{L4} (taking $p=2$) it is enough to show 
that $\{2,4, \ldots, 2k \}$ are in $\{a(1),\ldots, a(n) \}$.

Suppose not, and let $2m$ be the largest even number $\le 2k$ not in
$\{a(1), \ldots, a(n) \}$.
Every even number $a(i) > 2m$ with $i \le n$ will be followed by 
$a(i+1) \le 2k$ (since $2m$ is always available).
But in $\{a(1),\ldots, a(n) \}$ we have at least $4k$ even numbers, 
and so at least $2k$ even numbers $> 2k$.
Therefore in $\{a(1),\ldots, a(n+1) \}$ we see all the numbers from 1 to $2k$,
including $2m$, a contradiction.~~~$\bsq$

\begin{lemma}\label{L6}
If $a(n) =p$, a prime, then all numbers $\{1, \ldots, p-1\}$ 
occur in $\{a(1),\ldots, a(n-1)\}$.
\end{lemma}

\paragraph{Proof.}
By Lemma \ref{L1}, $a(n-1) = qp$, $q$ prime, $q< p$;
so $q$, $2q , \ldots, (p-1 ) q$ have already appeared in 
$\{a(1),\ldots, a(n-2) \}$.
The result now follows by Lemma \ref{L4}.~~~$\bsq$

\paragraph{Proof of Theorem~\ref{T3}.}
By direct verification we may assume $n \ge 260^2$.
Let $m = \lceil n/260 \rceil$, and
suppose, seeking a contradiction, that $a(n) < m$.
Note that the lower bound on $n$ implies that $n/m > 259$.

We will show that at least $4m$ even numbers have occurred in 
$\{a(1),\ldots, a(n-2) \}$, which gives a contradiction by
Lemma \ref{L5}.
No primes greater than $m$ can occur in this interval, 
or we get a contradiction by Lemma \ref{L6}.

Some number $\ge n$ must occur among $\{a(1),\ldots, a(n-1) \}$, 
since there are $n-1$ numbers and $a(n) < m$ is missing.
Suppose $a(j) \ge n$,
with controlling prime $p \le m$, say $a(j) = lp$, and $j \le n-1$.
Since $lp \ge n$, $l \ge 260$.
Let $R \le m$ be the smallest missing number among
$\{a(1), \ldots, a(n-1) \}$.
Then $l \le R$, for if $l > R$ then we have seen $Rp$ at time $j-1$,
and by Lemma \ref{L4} we have seen $R$ by time $j$, a contradiction.
Therefore $l \le m$ and so $p \ge n/l > 259$.
Both $l$ and $p$ are in the range $[260,m]$.

We consider the ``window''
$$\sW = \left\{ \left\lfloor \frac{lp}{4} \right\rfloor,
\left\lfloor \frac{lp}{4} \right\rfloor + p , \ldots, lp \right\} \,,
$$
and define ``entry'' and ``exit'' point as in the proof of Theorem \ref{T2}.

There are at least $\lfloor 3l/4 \rfloor$ multiples of $p$ in $\sW$, 
and at least
$\lfloor 3l/8 \rfloor -1$ of them are even.
Any such even multiple of $p$ is an exit point (from Lemma \ref{L5}).
There must therefore be at least
$\lfloor 3l/8 \rfloor -2$ entry points which are not controlled by $p$.
Let $\sS = \{q_1, q_2 , \ldots \}$ (with $q_1 < q_2 < \cdots \le l )$ 
be the set of controlling primes for these entry points.
Then
$$D_\sS \left( \left\lfloor \frac{l}{4} \right\rfloor, 
l \right) \ge \lfloor 3l/8 \rfloor -2 \,.
$$
As in \eqn{Eq8a}, \eqn{Eq8}, we get
\beql{Eq10}
\sum_{i \ge 1} \frac{1}{q_i} \ge
\frac{\left\lfloor \frac{3l}{8} \right\rfloor -2 - \pi (l)}
{l - \left\lfloor \frac{l}{4} \right\rfloor}
\eeq
which is
$\ge 43/214$ for $l \ge 260$.
We define $k$ by
\beql{Eq11}
\sum_{i=1}^{k-1} \frac{1}{q_i} < \frac{43}{214} \le 
\sum_{i=1}^k \frac{1}{q_i} ~.
\eeq

Every multiple of any element $q \in \sS$ that is 
$\le \lfloor lp/4 \rfloor$ must have already occurred in 
$\{a(1), \ldots, a(n-2) \}$.
Of these, at least $\lfloor lp / 8q \rfloor$ are even.
Therefore the number of distinct even numbers in the range 
$1, \ldots, \lfloor lp /4 \rfloor$ that have occurred is at least
$$\sum_{i \ge 1}^k \left\lfloor \frac{X}{2q_i}\right\rfloor -
\sum_{1 \le i < j \le k}
\left\lfloor \frac{X}{2q_i q_j}\right\rfloor
$$
where $X = lp /4$.
Proceeding as in the proof of Theorem \ref{T2}, 
and using $q_1 \ge 2$, we find that this is at least $4.024 m$.
This exceeds $4m$ and provides the desired contradiction.~~~$\bsq$

%
%

\section{Cycle structure}
Since the sequence is a permutation of $\NN$, it is also
natural to investigate the cycle structure.
The experimental evidence suggests
 that there are infinitely many finite cycles and 
infinitely many infinite cycles. It seems very likely to 
be hard to prove either of these observations, or even
to prove that there is at least one infinite cycle.

The first few finite cycles start at the points
$$1,2,3,8,40,64,121,149,359,2879,5563,28571,251677,$$
and have lengths
$$1,1,6,1,1,1,2,12,11,25,8,22,11$$
respectively.
There are also a large number of apparently 
infinite cycles, of which the first two are
$$\cdots 229310, 117833, \ldots, 22, 27, 26, 28, 13, 14, {\bf 7},
12, 18, 20, 11, 15, 21, \ldots, 636551, 652766, \ldots
$$
and
$$\ldots, 502008, 257519, \ldots, 248, 253, 131, 138, {\bf 73}, 82, 
129, 201, 212, \ldots, 645906, 662330, \ldots
$$
with minimal representatives 7 and 73 respectively.
The first fifteen of these apparently distinct
cycles have not coalesced in the 
first 700000 terms of the sequence.
However, although it seems unlikely, it is theoretically possible 
that they could coalesce at some later point.
It would be nice to know more!

%
%

\section{Generalizations}

The EKG sequence can be generalized in various ways,
while retaining the basic construction of a greedy sequence
with a condition on gcd's of consecutive terms. 
For fixed $M \ge 2$, let $b(n) = n$ for $1 \le n \le M$, and for 
$n \ge M+1$ let $b(n)$ be the smallest natural number not already 
in the sequence with the property that ${\rm gcd} \{b(n-1) , b(n) \} \ge M$.
The proof of Theorem \ref{T1} easily extends to show that 
$(b(n) : n \ge 1 )$ is also a permutation of $\NN$.
For the cases $M=3,4,5$ see  sequences \htmladdnormallink{A064417}{http://www.research.att.com/cgi-bin/access.cgi/as/njas/sequences/eisA.cgi?Anum=A064417}, \htmladdnormallink{A064418}{http://www.research.att.com/cgi-bin/access.cgi/as/njas/sequences/eisA.cgi?Anum=A064418}, \htmladdnormallink{A064419}{http://www.research.att.com/cgi-bin/access.cgi/as/njas/sequences/eisA.cgi?Anum=A064419}
in [Sloane 2001].

\begin{figure}[htb]
\centerline{\includegraphics[angle=270, 
width=5in]{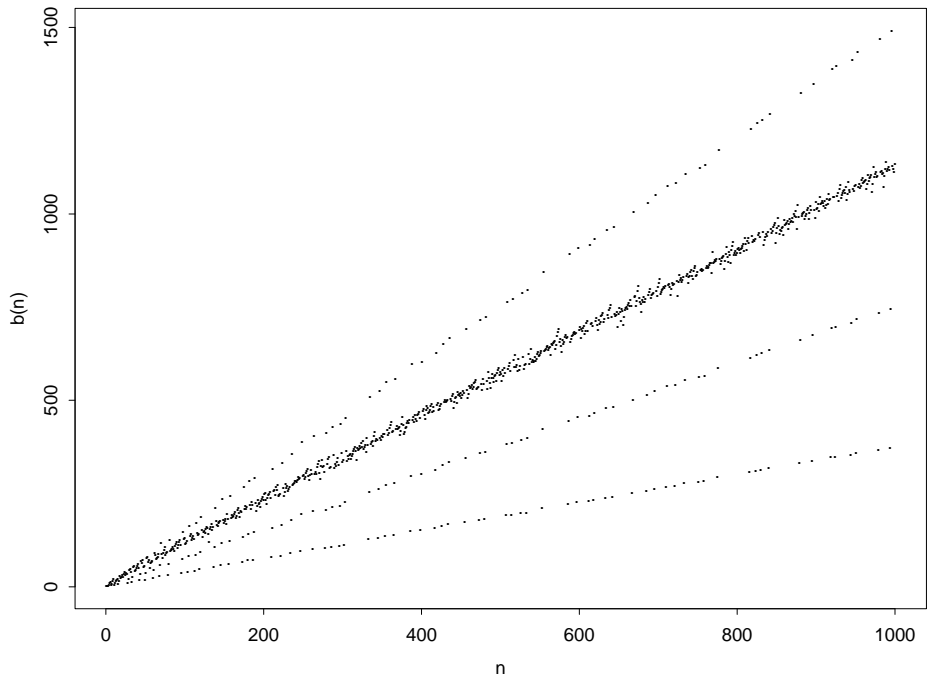}}

\caption{Plot of $b(1)$ to $b(1000)$ for the case $M=3$.}
\label{F7}
\end{figure}

The first 1000 terms for the case $M=3$ are shown in Figure~\ref{F7}.
This sequence appears to behave in a similar way to the EKG sequence: 
the spikes in the sequence are associated with
the primes, occurring in the order $3p, p, 2p, 4p$. All points of
the sequence seem to lie near lines of slope $1/3, 2/3, 1$ and $4/3$.

More generally, we can allow an arbitrary finite prefix at the beginning of
such sequences.
Consider a sequence $c(n)$ with the  generating rule
that ${\rm gcd} \{c(n-1) , c(n) \} \ge M$, for $n \ge N$, 
but with a finite
prefix $c(n)= a_n : 1 \le n \le N$, such that
all terms $a_n$ are distinct natural numbers, and
that the prefix includes values $a(n) = k$ for
of $1 \le k \le M$. With a little 
work the  proof of Theorem \ref{T1}
can be modified to apply to these sequences as well
to show they are permutations.

If we use the greedy gcd rule with $M=2$, starting with
an arbitrary  finite prefix, one obtains an infinite number 
of different EKG-like permutations. Many of these sequences,
perhaps all, eventually coalesce with the original EKG sequence.
In any case, the  qualitative properties of the
resulting permutations appear similar to the EKG sequence:
the primes provably appear in consecutive order (excluding those
in the prefix terms); the general plot of the permutations
remains similar to Figure~\ref{F3}.

It appears that linear upper and lower bounds hold for 
all the permutations $c(n)$ of these types.
Very likely such bounds can be established rigorously
in any particular case using methods similar to those used in
proving Theorems~\ref{T2} and \ref{T3}.

\subsection*{Acknowledgements}
We thank Jonathan Ayres for discovering this wonderful sequence.
We also thank a referee for helpful comments.
\section*{References}

\begin{description}
\item[{[Ayres 2001]}]
J. Ayres, personal communication, Sept. 30, 2001.

\item[{[Erd\H{o}s et~al. 1983]}]
P. Erd\H{o}s, R. Freud and N. Hegyvari,
Arithmetical properties of partitions of integers,
{\it Acta Math. Acad. Sci. Hungar.}, {\bf 41} (1983), 169--176.

\item[{[Hooley 1976]}]
C. Hooley,
{\em Applications of Sieve Methods to the Theory of Numbers},
Cambridge Tracts in Math. No. 70,
Cambridge Univ. Press: Cambridge 1976.

\item[{[Sloane 2001]}]
N. J. A. Sloane,
{\em The On-Line Encyclopedia of Integer Sequences},
published electronically at
\htmladdnormallink{www.research.att.com/$\sim$njas/sequences/}{http://www.research.att.com/~njas/sequences/}.

\item[{[Zagier 1977]}]
D. Zagier, The first 50 million prime numbers,
{\em Math. Intelligencer}, {\bf 0} (1977), 7--19.
\end{description}

\end{document}